\documentclass[letterpaper, 10 pt, conference]{ieeeconf}  

\IEEEoverridecommandlockouts                              
\usepackage{amsmath,amsgen,amstext,amssymb}
\usepackage{amsfonts}
\usepackage{dsfont}
\usepackage{mathrsfs}
\usepackage{latexsym}
\usepackage{times}
\usepackage{graphics}
\usepackage{graphicx}
\usepackage{epsfig}
\usepackage{psfrag}
\usepackage{bm}
\usepackage{lmodern}
\usepackage[T1]{fontenc}
\usepackage{url}
\usepackage{color}

\DeclareMathAlphabet{\mathpzc}{OT1}{pzc}{m}{it}
\DeclareMathOperator*{\argmax}{\arg\max}

\newtheorem{theorem}{Theorem}[section]

\newcommand{\beqn}{\begin{equation}}
\newcommand{\eeqn}{\end{equation}}

\newcommand{\ex}{{\mathbb{E}}}
\newcommand{\pr}{{\mathbb{P}}}

\newcommand\Reals{{\mathbb{R}}}
\newcommand\Ints{{\mathbb{Z}}}

\newcommand{\g}{\lambda}                

\newcommand{\cC}{{\mathcal{C}}}
\newcommand{\cD}{{\mathcal{D}}}
\newcommand{\cF}{{\mathcal{F}}}

\newcommand{\cP}{{\mathcal{P}}}

\newcommand{\cS}{{\mathcal{S}}}

\title{\LARGE \bf
Control of Time-Varying Epidemic-Like Stochastic Processes and Their Mean-Field Limits
}

\author{Yingdong Lu, Mark S.\ Squillante, Chai Wah Wu
\thanks{The authors are with the Mathematical Sciences Department at the IBM Thomas J.\ Watson Research Center, Yorktown Heights, NY 10598, U.S.A
        {\tt\footnotesize \{yingdong,mss,cwwu\}@us.ibm.com}}%
}

\begin{document}

 \maketitle
\thispagestyle{empty}
\pagestyle{empty}

\begin{abstract}
The optimal control of epidemic-like stochastic processes is important both historically and for emerging applications today,
where it can be especially important to include time-varying parameters that impact viral epidemic-like propagation.
We connect the control of such stochastic processes with time-varying behavior to the stochastic shortest path problem and obtain solutions for various cost functions.
Then, under a mean-field scaling, this general class of stochastic processes is shown to converge to a corresponding dynamical system.
We analogously establish that the optimal control of this class of processes converges to the optimal control of the limiting dynamical system.
Consequently, we study the optimal control of the dynamical system where the comparison of both controlled systems renders various important mathematical properties of interest.
\end{abstract}


%
\section{Introduction}
The mathematical analysis of epidemic-like behavior has a very rich and important history,
with the seminal work of Bernoulli on epidemiological models as a starting point in the
1766~\cite{Bern1766-simple}.
More recently, mathematical models of epidemic-like behavior have received considerable attention in the research literature,
which include models of various aspects of large networks such as the complex structures and behaviors of communication networks,
social media/networks, viral-propagation networks (e.g., epidemics, computer viruses and worms), and financial networks;
refer to, e.g.,
\cite{EasKle10}
and the references therein.
However, research on the control and optimization of such epidemic-like mathematical systems has been much more limited;
see, e.g., \cite{BoChGa+10}.
Even more importantly, this entire body of work has focused solely on static (non-time-varying) model parameters that impact
the complex structures and behaviors of the large epidemic-like systems of interest.
In contrast, time-varying behaviors often arise in many emerging applications of epidemic-like systems, especially those where one
observes behaviors that lead to forms of exacerbated complex dynamics and actions frequently found in communication, financial,
social, and viral-propagation networks.

We consider in this paper variants of the classical mathematical model of epidemic-like behavior analyzed by
Kurtz \cite[Chapter 11]{EthKur86},
extending the analysis to first incorporate time-varying behavior for the infection and cure rates of the model and to then
study aspects of the corresponding stochastic optimal control problem.
Specifically, we start by formally presenting an epidemic-like continuous-time, discrete-state stochastic process in which
each individual comprising the population can be either in a non-infected state or in an infected state, and where the rate at which
the non-infected population is infected and the rate at which the infected population is cured are both functions of time.
Then, we investigate the optimal control problem associated with this general class of epidemic-like stochastic processes.
Making connections to the well-studied stochastic shortest path problem, we exploit structural properties to obtain the corresponding
optimal control policy; in one special case of interest we derive the explicit control policy, whereas in other cases we compute
the value functions through efficient (linear program based) methods instead of value or policy iteration.
Taking the limit as the population size tends to infinity under a mean-field scaling, we establish that the time-varying epidemic-like
stochastic processes converge to a continuous-state nonautonomous dynamical system.
Moreover, the control of such stochastic processes is shown to be asymptotically equivalent to the optimal control of the limiting dynamical system.
We then investigate instances of the optimal control of the limiting dynamical system and establish structural properties of their equilibria and trajectories.
Lastly, computational experiments compare the optimal control of both the stochastic process as a function of population size and the dynamical system.

The paper is organized as follows.
\S\ref{sec:stochastic} presents our model and control of the general class of epidemic-like stochastic processes with time-varying parameters.
\S\ref{sec:convergence} presents results for the mean-field limit of the stochastic process and its controlled counterpart.
\S\ref{sec:limit} presents our model and analysis of the corresponding limiting dynamical system.
\S\ref{sec:exp} presents some computational experiments for both systems.
We refer the reader to \cite{LuSqWu17} for
additional results, all proofs, related work, and technical details.

\section{Epidemic-Like Stochastic Processes}
\label{sec:stochastic} 

\subsection{Mathematical Model}
Consider a sequence of Markov processes
$\hat{Z}_n = \{ (\hat{X}_n(t), \hat{Y}_n(t)) ; t \geq 0 \}$
indexed by the total population size $n \in \Ints^+ := \{1, 2, \ldots\}$ and defined over the probability space $(\hat{\Omega}_n, \cF_n, \pr_n)$,
composed of the state space $\hat{\Omega}_n := \{ (i,j) : 0 \leq i, j \leq n, i+j =n \}$, $\sigma$-algebra $\cF_n$ and
probability measure $\pr_n$, with initial probability distribution $\bm\alpha_n$.
Each process $\hat{Z}_n(t)$ represents the ordered pair $(\hat{X}_n(t),\hat{Y}_n(t))$ of non-infected and
infected population at time $t$, where we assume connections among the population form a complete graph.

Define $\Omega \subset \Reals^2$ and $\Omega_n := \Omega \cap \{ \ell/n : \ell \in \hat{\Omega}_n \}$.
The time-dependent infinitesimal generator $Q_n(t) = [q^{(n)}_{(i,j)(k,l)}(t)]_{ (i,j), (k,l) \in\hat{\Omega}_n}$ for the Markov process
$\hat{Z}_n$ has time-dependent transition intensities
$q^{(n)}_{(i,j)(i-1,j+1)}(t) = \lambda(t) i \frac{j}{n} = n \lambda(t) \frac{i}{n} \frac{j}{n}$,
$q^{(n)}_{(i,j)(i+1,j-1)}(t) = \mu(t) j = n \mu(t) \frac{j}{n}$,
where the latter equalities bear the general form $q_{k, k+\ell}^{(n)}(t) = n \beta_{\ell, t}(k/n)$, for $k, k+\ell \in \hat{\Omega}_n$,
with nonnegative functions $\beta_{\ell,t}(x)$ defined on $\Omega$ for $\ell \in \hat{\Omega}_n$ and $t \geq 0$, continuous in $t$,
Lipschitz continuous in $x=k/n$ (by definition), and $(x + \ell/n) \in \Omega_n$ when $\beta_{\ell, t}(x) > 0$, for $x \in \Omega_n$.
Throughout the functions $\lambda(t)$ and $\mu(t)$ are assumed to be continuous in $t$.

Note that the above definition of the epidemic-like stochastic process $\hat{Z}_n$ is slightly different from the corresponding
(non-time-varying) model of Kurtz~\cite{EthKur86}, in that we allow an infected individual who is cured to become infected at a later time.
In any case, our results hold for both types of epidemic-like models and related variations thereof with time-varying transition rates.

\subsection{Optimal Control and Analysis}
\label{sec:MDP:optimal}
Now consider the optimal control problem associated with the stochastic process $\hat{Z}_n$ for any population size $n$,
fixing $n$ henceforth in this subsection and omitting this parameter when clear by context.
Let $P(\hat{X}(t))$ denote the profits (rewards minus costs) as a function of the state of the system $\hat{X}(t)$ at time $t\in \Reals_+$.
The decision variables are based on the controlled infection and cure rates $\lambda(t)$ and $\mu(t)$ deployed by the system that represent
changes from the original infection and cure rates,
now denoted in this subsection by $\hat{\lambda}(t)$ and $\hat{\mu}(t)$,
where the system incurs additional costs $\hat{C}_\lambda(\lambda(t), \hat{\lambda}(t))$ and $\hat{C}_\mu(\mu(t), \hat{\mu}(t))$ as functions
of the pairs of infection and cure rates, respectively.
Throughout this subsection the control variables $\lambda(t)$ and $\mu(t)$ are assumed to be continuous in $t$,
with $\hat{\lambda}(t)$ and $\hat{\mu}(t)$ continuously varying for all $t$.
Define $\bm{\lambda} := (\lambda(t))$ and $\bm{\mu} := (\mu(t))$.
The objective function of our optimal control formulation is then given by
\begin{align}
\max_{\bm{\lambda}, \, \bm{\mu}} & \quad f\bigg( \, \int_{0}^T \Big\{ P(\hat{X}(t)) - \hat{C}_\lambda(\lambda(t), \hat{\lambda}(t)) \nonumber \\
& \qquad\qquad\qquad\qquad - \hat{C}_\mu(\mu(t), \hat{\mu}(t))\Big\} dt \, \bigg), \label{opt:obj}
\end{align}
where $T$ denotes the time horizon, which can be finite or infinite, and $f(\cdot)$ represents an operator of
interest.
An appropriate form of expectation is of primary interest in this paper.

Throughout this section, we assume that $P(x), \, x\in [0,c]$, has a single maximum at $x^*$ which, e.g., occurs when $P(\cdot)$
is linear (in which case $x^*=0$ or $x^*=c=1$) or when $P(\cdot)$ is concave (in which case $x^* \in [0,c=1]$).
Define $[n] := \{0,\ldots,n\}$ and let $i^*=\argmax_{i\in[n]} P(i)$ denote the integer(s) at which the profit function $P(\cdot)$
has a maximum value w.r.t.\ the states of $\hat{Z}$.
We further assume there is a unique $i^*$ in order to elucidate the exposition that follows;
otherwise, there would be $i^*_1 = \lfloor x^* \rfloor$ and $i^*_2 = \lceil x^* \rceil$ where $P(i^*_1) = P(i^*_2)$
and our analysis would apply w.r.t.\ both states $i^*_1$ and $i^*_2$.

To determine the optimal control policy for the Markov process $\hat{Z}$, we first exploit uniformization and consider the
uniform version discrete-time Markov chain $\hat{Z}^\prime$ with transition probability matrix $(p_{ij})_{i,j\in [n]}$
and uniformization constant $\nu > \max_{i \in [n]}\{-Q_{ii}\}$.
Our objective then is to determine the policy $\pi$ that maximizes the total expected profit over the entire time horizon
\begin{align*}
J_\pi(i) & = \ex\left[\sum_{t=1}^\infty h(i(t), u(i(t))) \Big|x_0=i\right], \\
J^*(i) & = \max_\pi J_\pi(i) ,
\end{align*}
or equivalently $J^*(i) = \min_\pi J_\pi(i)$ and
\begin{equation}
J_\pi(i) = \ex\left[\sum_{t=1}^\infty g(i(t), u(i(t))) \Big|x_0=i\right],
\label{eqn:genBellman}
\end{equation}
where $-h(i(t),u(i(t)))=g(i(t),u(i(t)))=-P(i(t)) +\hat{C}_\lambda(\lambda(i(t)),\hat{\lambda}(t)) +\hat{C}_\mu(\mu(i(t)),\hat{\mu}(t))$
denotes the cost per stage for state $i(t)$ of $\hat{Z}^\prime$ at time $t\in \Ints^+$, for all $i\in[n]$.

Next, the above expressions and assumptions allow us to determine the optimal control policy by solving the optimal control problem
through a corresponding stochastic shortest path problem~\cite{Bertsekas} over the set of states $\cS \cup \{ i^* \}$ where
$\cS := [n] \setminus \{ i^* \}$ and state $i^*$ is the special cost-free absorbing state in which the system remains at no
further cost once reached.
For all other states $i\in \cS$, a cost of $g(i,u)$ will be incurred if action $u\in U(i)$ is taken when the Markov chain $\hat{Z}^\prime(t)$
is in state $i$ at time $t\in \Ints^+$.
The objective is to determine the policy $\pi$ such that, for all $i\in \cS$, $J^*(i) = \min_\pi J_\pi(i)$ together with \eqref{eqn:genBellman}.
From known results for the stochastic shortest path problem, specifically Proposition 7.2.1 in \cite{Bertsekas}, we know that the optimal costs
$J^*(i), \, \forall i \in \cS$, satisfy Bellman's equation
\begin{equation}
J^*(i) = \min_{u\in U(i)} \left[g(i,u) + \sum_{j \in [n]} p_{ij}(u) J^*(j)\right]
\label{eqn:Bellman}
\end{equation}
and that a stationary policy $\gamma$ is optimal if and only if for every state $i\in \cS$, $\gamma(i)$ attains the minimum in \eqref{eqn:Bellman}.

While the value functions and optimal policies in general can be computed through methods such as value or policy iteration, the remainder of this
section considers several cases in which exact solutions can be derived or exact calculations can be applied to efficiently obtain the value functions.
Our analysis will exploit structural properties of the Markov chain $\hat{Z}^\prime$, such as
$p_{ij}=0$ for all $j \notin \{i-1,i,i+1\}$ when $i\in\{1,\ldots,n-1\}$;
$p_{0j}=0$ for all $j \notin \{0,1\}$; and
$p_{nj}=0$ for all $j \notin \{n-1,n\}$.

\subsubsection{No Action Costs}
We start with the case where there are no costs for adjusting the infection and cure rates, i.e., $\hat{C}_\lambda(a,b) = 0 = \hat{C}_\mu(a,b)$ for all $a,b$.
Hence, $g(i,u) = -P(i) = c(i)$ for any $(i,u)$, $i\in \cS, \, u \in U(i)$, and thus $c(i)$ is monotone decreasing when $i < i^*$ and monotone increasing
when $i > i^*$.
In conjunction with the special structure of the Markov chain, Bellman's equation for this case can be written as 
\begin{align}
J(i) = c(i) &+ \min_{\g\in[0, {\bar \g}], \mu\in [ 0, {\bar \mu}] }\left[\frac{\mu (n-i)}{\nu}J(i-1) \right. \nonumber \\
& \left. +\left(1-\frac{\mu (n-i)}{\nu}-\frac{\g i(n-i)}{n\nu}\right)J(i) \right. \nonumber \\
&\left. + \frac{\g i(n-i)}{n\nu}J(i+1)\right]. \label{eqn:OurBellman}
\end{align}

From the properties of $c(i)$ and an analysis of the minimization in Bellman's equation, we observe that the optimal policy should push the Markov chain
$\hat{Z}^\prime$ to the right as hard as possible when $i < i^*$ and push to the left as hard as possible when $i > i^*$.
Based on this observation, we can ``guess'' an optimal control policy and the corresponding set of $J(i), \, \forall i\in\cS$, and then verify that this
system of equations satisfies \eqref{eqn:OurBellman}.

In particular, we respectively have for $i \leq i^*$ and $i \geq i^*$
\begin{align*}
J(i) & = c(i) + \left[\frac{\mu (n-i)}{\nu}J(i-1) +\left(1-\frac{\mu (n-i)}{\nu}\right)J(i) \right], \\
J(i) & = c(i) + \left[ \left(1-\frac{\g i(n-i)}{n\nu}\right)J(i) + \frac{\g i(n-i)}{n\nu}J(i+1)\right],
\end{align*}
where in both systems of equations we set one of the control rates to be zero and set one of the other control rates to be its maximum.
We therefore obtain
\begin{align}
J(i) & = J(i+1) + \frac{n\nu c(i)}{\g i(n-i)}, \quad i < i^*,  \label{eqn:recursion_1} \\
J(i) & = J(i-1) + \frac{\nu c(i)}{\mu(n-i)}, \quad i > i^*.  \label{eqn:recursion_2}
\end{align}
It then can be readily verified that \eqref{eqn:recursion_1} and \eqref{eqn:recursion_2} are indeed a solution to Bellman's equation \eqref{eqn:OurBellman},
which from Proposition 7.2.1 in \cite{Bertsekas} is the unique solution of \eqref{eqn:OurBellman}.

\subsubsection{Linear Action Costs}
Let us now consider the case where the costs for adjusting the infection and cure rates are linear.
Specifically, to simplify the presentation, we assume the action costs to be linear functions of $\lambda$ and $\mu$, i.e.,
$\hat{C}_\lambda(\lambda(t), \hat{\lambda}(t)) = c_\lambda \cdot \lambda(t)$ and $\hat{C}_\mu(\mu(t), \hat{\mu}(t)) = c_\mu \cdot \mu(t)$
where $c_\g$ and $c_\mu$ are two constants for the corresponding cost rates.
In this case, Bellman's equation becomes
\begin{align*}
J(i) = c(i) &+ \min_{\g\in[0, {\bar \g}], \mu\in [ 0, {\bar \mu}]}\left[c_\g\g + c_\mu \mu+  \frac{\mu (n-i)}{\nu}J(i-1) \right. \\
&+\left(1-\frac{\mu (n-i)}{\nu}-\frac{\g i(n-i)}{n\nu}\right)J(i) \\
&+\left. \frac{\g i(n-i)}{n\nu}J(i+1)\right].
\end{align*}

We can rewrite this expression to obtain
\begin{align*}
0 = c(i) &+ \min_{\g\in[0, {\bar \g}], \mu\in [ 0, {\bar \mu}]}\left[\left(c_\g+\frac{i (n-i)}{n\nu}\Delta(i+1)\right) \g \right. \\
&+ \left. \left(c_\mu - \frac{(n-i)}{\nu} \Delta(i)\right)\mu \right] ,
\end{align*}
where $\Delta(i) := J(i) -J(i-1)$.
This optimization problem is clearly a linear program, which implies that only the vertices of the feasible region need to be considered.
The vertices consist of the $\g$ and $\mu$ that take on values of either $0$ or its maximum ${\bar \g}$ and ${\bar \mu}$, respectively.
We therefore can further reduce Bellman's equation as follows
\begin{align*}
0 = c(i) &+ \min_{\g\in\{0, {\bar \g}\}, \mu\in \{ 0, {\bar \mu}\}} \left[\left(c_\g+\frac{i (n-i)}{n\nu}\Delta(i+1)\right) \g \right. \\
&+ \left. \left(c_\mu - \frac{(n-i)}{\nu} \Delta(i)\right)\mu \right].
\end{align*}
This represents Bellman's equation for a stochastic control problem with finite states and finite controls,
for which it is well known~\cite{Bertsekas} that the problem can be solved as a linear program. 
Hence, we efficiently compute the solution of our control problem in this case via 
\begin{align*}
\max \quad & \sum_i \eta_i J(i)  \\
s.t.  \quad & 0 \le c(i) +  \left[\left(c_\g+\frac{i (n-i)}{n\nu}(J(i+1)-J(i))\right) \g \right. \\
\quad & \qquad \qquad \left.+ \left(c_\mu - \frac{(n-i)}{\nu} (J(i)-J(i-1))\right)\mu \right], \\
\quad & \qquad \qquad \qquad \forall \g\in\{0, {\bar \g}\}, \mu\in \{ 0, {\bar \mu}\} .
\end{align*}
where $\eta_i$ are positive real numbers and $J(i)$ are the variables for the linear program.

\subsubsection{General Action Costs}
Lastly, consider the case where the infection and cure rate adjustment costs are general functions with the number of actions
restricted to a finite number of possibilities, in which case we can extend our above linear program approach.
More specifically, Bellman's equation in this case can be written as
\begin{align}
0 = c(i) &+ \min_{(\g, \mu) \in U_K} \left[\left(c_\g+\frac{i (n-i)}{n\nu}\Delta(i+1)\right) \g \right. \nonumber \\
&+ \left. \left(c_\mu - \frac{(n-i)}{\nu} \Delta(j)\right)\mu \right], \label{eqn:OurBellman_linear_discrete}
\end{align}
where $U_K$ represents the finite set of $K$ possible combinations of $\g$ and $\mu$.
The value functions then can be obtained by solving the following linear program 
\begin{align*}
\max \quad & \sum_i \eta_i J(i)  \\
s.t.  \quad & 0 \le c(i) +  \left[\left(c_\g+\frac{i (n-i)}{n\nu}\Delta(i+1)\right) \g \right. \\
\quad & \quad \left.+ \left(c_\mu - \frac{(n-i)}{\nu} \Delta(j)\right)\mu \right], \quad \forall (\g,\mu)\in U_K .
\end{align*}

In general, we can always discretize the action space to obtain a version of \eqref{eqn:OurBellman_linear_discrete},
the solution of which serves as an approximation to the value functions of the original control problem whose accuracy
can increase with $K$.

\section{Mean-Field Limits}
\label{sec:convergence}

\subsection{Epidemic-Like Stochastic Processes}
Suppose that the Markov Chain $\hat{Z}_n(t)$ is as defined in \S\ref{sec:stochastic} with time-dependent transition intensities of the general form given therein.
From the martingale-problem method (see, e.g., \cite[Chapters~4,~6]{EthKur86}), we devise that $\hat{Z}_n(t)$ has the integral representation
\begin{align} \label{eqn:int_rep}
\hat{Z}_n (t) & = \hat{Z}_n(0) + \sum_\ell \ell W_\ell \left( n \int_0^t  \beta_{\ell, s}\left(\frac{\hat{Z}_n(s)}{n}\right) ds\right),
\end{align}
where the $W_\ell$ are independent standard Poisson processes.
Define $F_t (z) := \sum_\ell \ell \beta_{\ell, t}(z)$, $z \in \Omega_n$.
Further define $Z_n(t) := \hat{Z}_n(t)/n$ on the state space $\Omega_n$ with
transition intensities
$q_{i,j}^{(n)}(t) = n\beta_{n(j-i),t}(i)$, $i,j \in \Omega_n$.

Our strategy for the
proof is to first obtain the integral representation of $Z_n(t)$, which leads to the generator of $Z_n(t)$
again through the martingale-problem method and the law of large numbers
for the Poisson process.
From this and the above we derive the desired expression
\begin{align} \label{eqn:scaled_int_rep}
Z_n(t) \; = \; Z_n(0) &+ \sum_\ell  \frac{\ell}{n} \bar{W}_\ell \left( n \int_0^t \beta_{\ell, s} (Z_n(s)) ds \right) \nonumber \\
&+ \int_0^t F_s(Z_n(s))ds,
\end{align}
where $\bar{W}_\ell$ denotes the centered Poisson process, i.e., $\bar{W}_\ell(x)= W_\ell (x)-x$.
%
One of our main results can now be presented, upon noting the following basic fact:
$\lim_{n\rightarrow \infty} \sup_{u\le v} \|\frac{\bar{W}_\ell(nu) }{n} \| =  0$, a.s. for $v\ge 0$.

\begin{theorem}
        Suppose that for each compact set $K \subset \Omega$
$\sum_\ell |\ell| \sup_{x\in K} \beta_{\ell,t} (x) < \infty$, $\forall t\ge 0$,
        and there exists $M_K>0$ s.t.
        \begin{equation}
        \label{eqn:lip}
        |F_t(x)-F_t(y) |\le M_K|x-y|, \qquad \forall x,y \in K, t \ge 0.
        \end{equation}
        Further supposing $Z_n(t)$ satisfies \eqref{eqn:scaled_int_rep} and $\lim_{n\rightarrow \infty} Z_n(0) = z_0$, and denoting $Z(t)$ as the solution to
        \begin{equation}
                Z(t) \; = \; z_0 + \int_0^t F_s(Z(s)) ds, \qquad t\ge 0,
        \label{eq:Z-dynamics}
        \end{equation}
        then we have, for every $t\ge 0$,
        \begin{equation}
        \label{eqn:main_conv}
        \lim_{n\rightarrow \infty}\sup_{s\le t} |Z_n(s) -Z(s)| =0, \qquad a.s.
        \end{equation}
        \label{thm:Kurtz-new2}
\end{theorem}

From Theorem~\ref{thm:Kurtz-new2}, we then have that the stochastic process $Z_n(t)$ converges to a deterministic
process $Z(t)$, taking values in $[0,c]$, a.s.\ as $n\rightarrow \infty$ and that $Z(t)$ satisfies the integral form of the general nonautonomous dynamical system given in
\eqref{eq:Z-dynamics} where the specific details of the process and the corresponding set of ODEs depend upon $F_s(\cdot)$ characterizing the averaging behavior of the original
stochastic process $\hat{Z}_n(t)$.
For epidemic-like models, the process $Z_n(t)$ converges to a deterministic process $Z(t)=(X(t), Y(t))$ a.s.\ as $n\rightarrow \infty$
with $Z(t)$ satisfying the pair of ordinary differential equations (ODEs):
$\frac{dX(t)}{dt} = -\lambda(t) X(t) Y(t) + \mu(t) Y(t)$,
$\frac{dY(t)}{dt} = \lambda(t) X(t) Y(t) - \mu(t) Y(t)$, with proper initial conditions.
This desired a.s.\ convergence result justifies the use of a continuous-state nonautonomous dynamical system to model a discrete-state real-world stochastic system.

\subsection{Controlled Stochastic Processes}
We next turn our attention to an optimal control problem associated with the class of epidemic-like stochastic processes, where our goal is to show that this
control process is asymptotically equivalent to the optimal control of the corresponding set of ODEs as the population size tends to infinity under a mean-field scaling.

Consider a sequence of general controlled Markov processes $\hat{Z}_n(t)$, with the adaptive control process $u_n(t)$ that is realized w.r.t.\ the adaptive transition kernel
$n \beta_{\ell,t}(k/n)$, $k,k+\ell \in \hat{\Omega}_n$, recalling $\beta_{\ell,t}(\cdot)$ is continuous in $t$.
For each system indexed by $n$, the optimal control $u^*_n(t)$ is determined by solving the optimal control problem w.r.t.\ the cost functions $c_1(\cdot)$ and $c_2(\cdot)$:
\begin{align*}
	\hat{J}^*_n(z) = &\min_{u_n(t)} \quad \hat{J}_n(z) \\
	= &\min_{u_n(t)} \quad \left\{ \int_0^T  c_1(\hat{Z}_n(t), u_n(t) ) dt + c_2(\hat{Z}_n(T))\right\}, \\
	&\mbox{ s.t. } \quad \hat{Z}_n(0) = z.  
\end{align*}
Here we assume the cost functions $c_1(z,u)$ and $c_2(z)$ are uniformly bounded, which is reasonable and justified by our interest in costs related only to the proportion of a population.
Recall the integral representation of $\hat{Z}_n(t)$ and $Z_n(t)$ in \eqref{eqn:int_rep} and \eqref{eqn:scaled_int_rep}, respectively.

For comparison towards our goal in this section, we also consider the corresponding optimal control problem associated with the limiting mean-field dynamical system of the above subsection.
Namely, the optimal control $u^*(t)$ is determined by solving the corresponding optimal control problem w.r.t.\ the same cost functions $c_1(\cdot)$ and $c_2(\cdot)$,
which can be formulated as
\begin{align*}
	J^*(z) & = \min_{u_n(t)} \quad J(z) \\
	& = \min_{u_n(t)} \quad \left\{ \int_0^T  c_1(Z(t), u(t) ) dt + c_2(Z(T))\right\} , \\
	\mbox{ s.t. } & \quad Z(0) = z ,
\end{align*}
where $Z(t)$ follows the dynamics
$Z(t) = z + \int_0^t F_s(Z(s)) ds$.
Note that the function $F_s(\cdot)$ encodes the control information. 

We seek to show that the optimal control $u^*(t)$ in the limiting mean-field dynamical system provides an asymptotically equivalent optimal control $u^*_n(t)$
for the original system indexed by $n$ in the limit as $n$ tends toward infinity.
More specifically, we first establish the following main result.
\begin{theorem}\label{thm:control1}
	Let $\hat{Z}_n(t)$, $Z_n(t)$ and $Z(t)$ be as above.
	We then have
	\begin{align}
		\label{eqn:main}
		\lim_{n\rightarrow \infty} \hat{J}^*_n(z) = J^*(z).
	\end{align}
Furthermore, let $F_s^*(\cdot)$ denote the function that encodes the optimal control $u^*(t)$ of the limiting mean-field dynamical system.
Suppose the original stochastic process $\hat{Z}_n(t)$ follows the deterministic state-dependent control policy determined by $F_s^*(\cdot)$.
Then, asymptotically as $n\rightarrow\infty$ under a mean-field scaling, both systems will realize the same objective function value in \eqref{eqn:main}.
\end{theorem}

\section{Epidemic-Like Dynamical Systems}
\label{sec:limit} 

\subsection{Mathematical Model and Analysis}
We next consider the continuous-time, continuous-state nonautonomous dynamical system $z(t)=(x(t),y(t))$ from the results of Section~\ref{sec:stochastic}.
The starting state $(x(0),y(0))$ of the system at time $t=0$ has initial probability distribution $\bm\alpha$.
We assume throughout that $\lambda(t) > 0$ and $\mu(t) > 0$.
The state equations are then given by
$\frac{dx(t)}{dt}= -\lambda(t) x(t) y(t) + \mu(t) y(t)$ and $\frac{dy(t)}{dt}=  \lambda(t) x(t) y(t) - \mu(t) y(t)$,
where $x(t)$ and $y(t)$ respectively describe the fraction of non-infected and infected population at time $t$, 
with total population $c = x(t)+y(t)$.

Since $c = x(t)+y(t)$ and $\frac{d(x(t)+y(t))}{dt} = 0$, we have $x(t)+y(t)= c = x(0)+y(0)$ for all $t$; i.e., the total population is constant.
Upon substituting $y(t) = c-x(t)$, we can equivalently rewrite the two-dimensional ODE as an one-dimensional ODE:
$\frac{dx(t)}{dt} = \lambda(t) x(t)^2 - (\lambda(t) c+\mu(t))x(t) + \mu(t) c$.
We then have the following main result.
\begin{theorem}
        For the dynamical system $(x(t),y(t))$ with $0 \leq x(0), y(0) \leq c$ and $x(t)+y(t)= c$, $\lambda(t)$, $\mu(t)$ continuously varying for all $t$,
        the system has an asymptotic state at $x_1^*(t)=\frac{\mu(t)}{\lambda(t)}$ and an equilibrium point at $x_2^*=c$,
and stability properties given as follows.
        \begin{enumerate}
                \item $0 < \frac{\mu(t)}{\lambda(t)} < \xi < c$, $\forall t$: The equilibrium point $x_2^*$ is unstable.
                All trajectories of the dynamical system with initial state $x(0) < c$ will converge towards being
                eventually near the asymptotic state $x_1^*(t)$ w.r.t.\ a $\delta$-neighborhood, i.e.,
                $\|x(t) - \frac{\mu(t)}{\lambda(t)}\| \leq \delta$
                where $\delta$ is a nonnegative constant that depends on the rates of change of $\mu(t)$ and $\lambda(t)$.
                \item $\frac{\mu(t)}{\lambda(t)} > c$: The equilibrium point $x_2^*$ is stable, towards which
                all trajectories of the dynamical system will converge.
                \item $\frac{\mu(t)}{\lambda(t)} = c$: There is one equilibrium point at $x_2^*$, which is neither stable nor unstable,
                towards which all trajectories of the dynamical system will converge.
                \item $\frac{\mu(t)}{\lambda(t)} = 0$: There is one equilibrium point at $x_1^*=0$, which is neither stable nor unstable,
                towards which all trajectories of the dynamical system will converge.
        \end{enumerate}
        \label{thm:CTCS-vary}
\end{theorem}

To summarize, for the dynamical system of Theorem~\ref{thm:CTCS-vary}, all trajectories $x(t)$ will approach a $\delta$-neighborhood of
$\frac{\mu(t)}{\lambda(t)}$ when $0 < \frac{\mu(t)}{\lambda(t)} < c$, will approach $0$ when $\frac{\mu(t)}{\lambda(t)}=0$,
and will approach $c$ when $\frac{\mu(t)}{\lambda(t)} \geq c$.

\subsection{Optimal Control and Analysis}
\label{sec:DS:optimal}
Now we turn to consider the optimal control problem formulation \eqref{opt:obj} within the context of the dynamical system $z(t)$,
recalling our use of $\lambda, \hat{\lambda}, \mu, \hat{\mu}$ from \S\ref{sec:MDP:optimal}.
Let $\bm{\lambda}^*$ and $\bm{\mu}^*$ denote the optimal solution to \eqref{opt:obj} subject to the corresponding ODEs of the above subsection.
This formulation represents the general case of the optimal control problem of interest for the dynamical system.
Although there are no explicit solutions in general, this problem can be efficiently solved numerically using known methods from control theory.

To consider more tractable cases, and gain fundamental insights into the problem,
we start by first considering a one-sided version of this general problem in equilibrium with a fixed constant infection rate $\lambda=\hat{\lambda}=\hat{\lambda}(t)$
where the goal is to maximize the reward at the equilibrium point and only the parameter $\mu$ is under our control.
The optimal control in this case is a stationary policy for the cure rate, i.e., a single control $\mu$ in equilibrium.
Under a linear profit function with rate $\cP$ and linear action-cost function of $\mu$ with rate $\hat{\cC}_\mu$,
we can rewrite the objective function \eqref{opt:obj} as
\begin{equation*}
\max_{\mu} \;\; \cP(x(\infty)) - \hat{\cC}_\mu(\mu) ,
\end{equation*}
since the optimal control is a stationary policy for the cure rate.
Upon substituting
$\min\{ c , \frac{\mu}{\lambda} \}$ for $x(\infty)$,
we derive the optimal control policy to be
\begin{align}
\mu^* = \argmax_{\mu \geq 0} \; \cP\bigg(\min\Big\{ c , \frac{\mu}{\lambda} \Big\} \bigg)
- \hat{\cC}_\mu(\mu) . \label{opt:one-sided:solution}
\end{align}
Namely, the optimal stationary control policy employs for all time $t$ the single control $\mu^*$ that solves \eqref{opt:one-sided:solution}.
An analogous formulation and result on $\lambda^*$ can be established for the opposite one-sided version of the problem in
equilibrium with constant cure rate $\mu$.

Next, as another step toward the general formulation,
consider the case where there are no costs for adjusting the infection and cure rates, i.e., $\hat{C}_\lambda(a,b) = 0 = \hat{C}_\mu(a,b)$ for all $a,b$.
Further assume, as in \S\ref{sec:MDP:optimal}, that $P(x)$ continues to have a single maximum at $x^*$.
We introduce the notion of an {\em ideal trajectory} denoted by $(x^I(t) = x^*)$
that maximizes the objective function~(\ref{opt:obj}) at all time in this problem instance.
Hence, the optimal policy is to have $\frac{\mu(t)}{\lambda(t)} = x^*$ with $\lambda(t)$ as large as possible,
subject to $\frac{\hat{\mu}(t)}{\hat{\lambda}(t)}$ varying over time, since this governs the speed at which
$x(t)$ approaches and continually follows $x^*$.

%
More precisely, we now present a main result of interest for this instance of the general formulation showing that we can get arbitrarily close to the ideal trajectory, and thus the maximum objective.
\begin{theorem}
        Suppose $\hat{C}_\lambda(a,b) = 0 = \hat{C}_\mu(a,b)$, for all $a,b$.
        For each $\epsilon > 0$ with $0\leq x(0) < c-\epsilon$, there is a $\hat{\delta} > 0$ s.t.\ if $\lambda(t), \mu(t) > \hat{\delta}$
        and $\frac{\mu(t)}{\lambda(t)} = x^*$ for all $t$, then the optimal solution of \eqref{opt:obj} is realized within $\epsilon$.
        \label{thm:optimal}
\end{theorem}

Let us now consider the above case where there are no costs for adjusting the infection and cure rates,
but where there are constraints on the rates of change of the control variables $\lambda(t)$ and $\mu(t)$,
i.e., $\theta_\lambda^\ell < \dot{\lambda} < \theta_\lambda^u$ and $\theta_\mu^\ell < \dot{\mu} < \theta_\mu^u$.
We again assume that $P(x)$ has a single maximum at $x^*$.
Our above notion of an {\em ideal trajectory} remains the same,
namely $(x^I(t) = x^*)$ maximizes the objective function~(\ref{opt:obj}) without constraints for all time $t$.
We therefore have that the optimal policy consists of setting $\lambda(t)$ and $\mu(t)$ so as to maximize the speed at which
$x(t)$ approaches and continually follows a maximum within an achievable neighborhood of $x^*$, subject to the constraints on
$\dot{\lambda}$ and $\dot{\mu}$ and subject to $\frac{\hat{\mu}(t)}{\hat{\lambda}(t)}$ varying over time.

More precisely, we establish a result showing that we can get arbitrarily close to a best state within a $\delta$-neighborhood of the ideal trajectory,
and thus the maximum objective, where $\delta$ is a nonnegative constant that depends on the rates of change of $\hat{\lambda}(t)$ and $\hat{\mu}(t)$,
and on $\theta_\lambda^\ell , \theta_\lambda^u, \theta_\mu^\ell , \theta_\mu^u$.
Define $\bar{\cD} := \{ x(t) : \|x(t) - x^*\| \leq \delta \}$ and $\cD(t) := \{ x(t) : x(t) \in \bar{\cD} \mbox{ \it{and} } x(t) \mbox{ \it{is reachable at time} } t \}$ for all $t$.
The main result of interest for this instance of the general formulation can then be expressed as follows.
\begin{theorem}
        Suppose $\hat{C}_\lambda(a,b) = 0 = \hat{C}_\mu(a,b)$, for all $a,b$, together with the constraints
        $\theta_\lambda^\ell < \dot{\lambda} < \theta_\lambda^u$ and $\theta_\mu^\ell < \dot{\mu} < \theta_\mu^u$.
        For each $\epsilon > 0$ with $0\leq x(0) < c-\epsilon$, there is a $\hat{\delta} > 0$ s.t.\ if $\lambda(t), \mu(t) > \hat{\delta}$
        and $\frac{\mu(t)}{\lambda(t)} = \hat{x}^*(t) := \argmax_{x(t) \in \cD(t)} P(x(t))$ for all $t$,
        then the optimal solution of \eqref{opt:obj} under the constraints on $\dot{\lambda}$ and $\dot{\mu}$ is realized within
        a $\delta$-neighborhood of $x^*$, i.e., $\bar{\cD}$, and in particular, the optimal reachable solution $\hat{x}^*(t)$
        is realized within $\epsilon$.
        \label{thm:optimal2}
\end{theorem}

When the costs for adjusting the infection and cure rates are introduced to either of the above instances of the general formulation,
the optimal policy will deviate from the ideal policies above where the deviation will depend on the initial state $x(0)$,
the cost functions $\hat{C}_\lambda(\cdot,\cdot)$ and $\hat{C}_\mu(\cdot,\cdot)$, the rates of change of $\hat{\lambda}(t)$ and $\hat{\mu}(t)$,
and any constraints on the rates of change of $\lambda(t)$ and $\mu(t)$.
Even though the policy of following the ideal trajectory is not optimal in general, it can provide structural properties and insight
into the complex dynamics of the system in a very simple and intuitive manner.

\section{Computational Experiments}
\label{sec:exp}
In this section we investigate various aspects of our theoretical results through computational experments.
The behavior of the Markov decision process is clear from the results in \S\ref{sec:MDP:optimal};
similarly for the behavior of the optimal control of the dynamical system from the results in \S\ref{sec:DS:optimal}.
We assume that the profit function $P$ is continuous and has a single maximum at $x^*$.
Then, when sampled in $n$ equal-spaced intervals, the set of $P(i/n)$ for $i=1,\cdots ,n$ reaches its maximum at either one or two values of $i$.
For simplicity, let us assume that for each $n$, $P(i/n)$ is maximized at a single value of $i^*$ such that $i^*/n$ is closest to $x^*$.
Classical number theory shows that the difference between $i^*/n$ and $x^*$ is asymptotically no better than $o(n^{-2})$ and this optimal rate is approached by the convergents of $x^*$.
Hence, the Markov decision process will converge towards this value of $i^*$.

\begin{figure}[htbp]
\centerline{\includegraphics[width=3.4in]{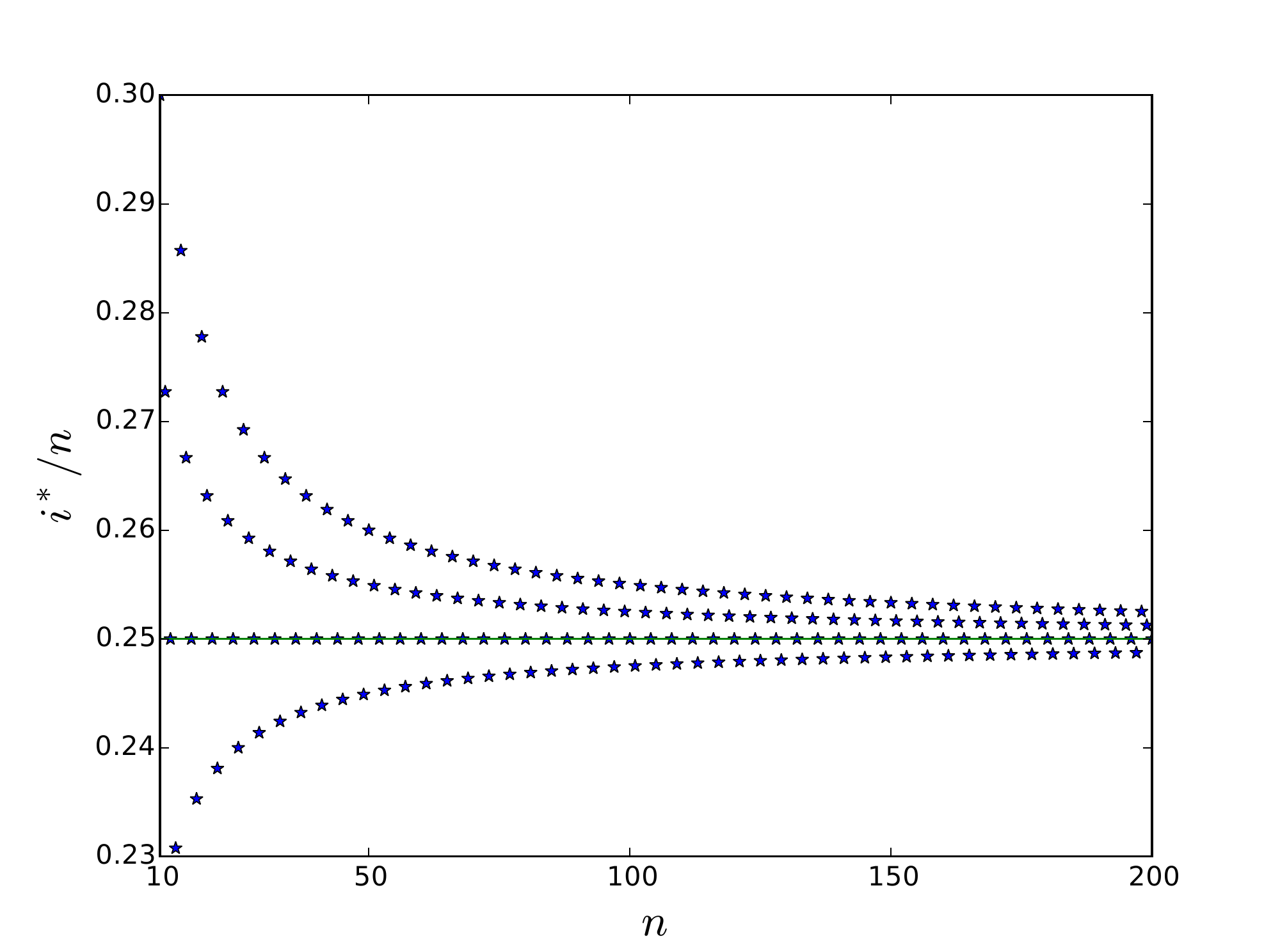}}
\caption{$i^*/n$ versus $n$, where $i^*$ is the state the MDP converges to. $x^* = \frac{1}{4}$ is indicated by the solid line.}\label{fig:figure1}
\end{figure}

\begin{figure}[htbp]
\centerline{\includegraphics[width=3.4in]{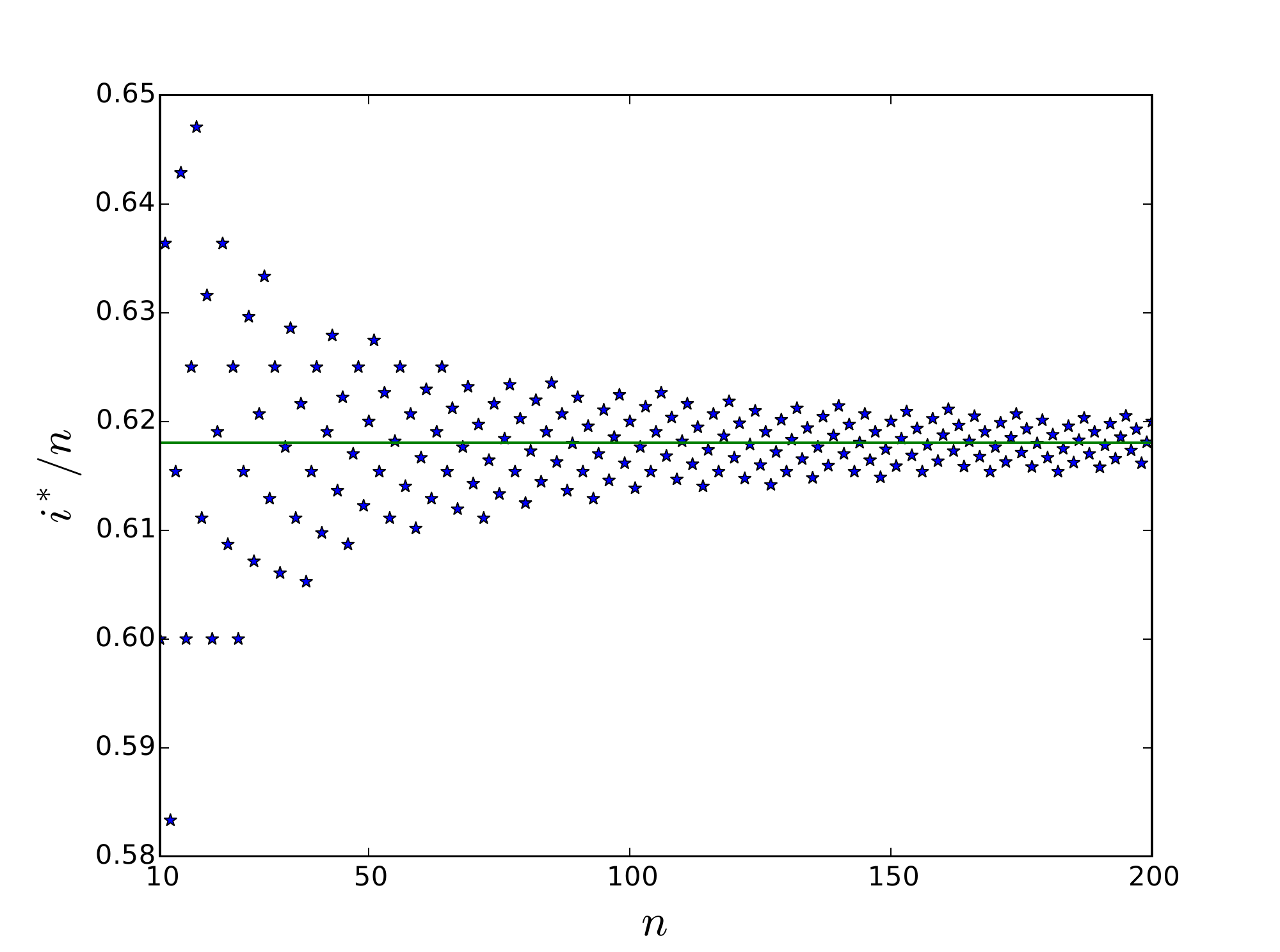}}
\caption{$i^*/n$ versus $n$, where $i^*$ is the state the MDP converges to. $x^* = \phi -1$ is indicated by the solid line.}\label{fig:figure2}
\end{figure}

Figures \ref{fig:figure1} and \ref{fig:figure2} illustrate how the value of $i^*$ behaves for $x^* = 1/4$ and $x^* = \phi - 1 = \frac{\sqrt{5}-1}{2}$, respectively.
We observe that the quantitative differences between the optimal control of the stochastic process and the dynamical system vanishes as $n \rightarrow \infty$,
and does so relatively quickly in accordance with classical results.

%

\begin{thebibliography}{10}



\bibitem{Bern1766-simple}
D.~Bernoulli.
\newblock Essai d'une nouvelle analyse de la mortalite causee par la
  petite verole.
\newblock {\em Mem.\ Math.\ Phys.\ Acad.\ Roy.\ Sci., Paris}, pages 1--45,
  1766.

\bibitem{Bertsekas}
D. ~Bertsikas
\newblock {\em Dynamic Programming and Optimal Control, Vols.\ I and II},
\newblock Athena Scientific, 2005.

\bibitem{BoChGa+10}
C.~Borgs, J.~Chayes, A.~Ganesh, A.~Saberi.
\newblock How to distribute antidote to control epidemics.
\newblock {\em Random Structures \& Algorithms}, 2010.



\bibitem{EasKle10}
D.~Easley, J.~Kleinberg.
\newblock {\em Networks, Crowds, and Markets: Reasoning About a Highly
  Connected World}.
\newblock Cambridge University Press, 2010.

\bibitem{EthKur86}
S.~N. Ethier, T.~G. Kurtz.
\newblock {\em Markov Processes: Characterization and Convergence}.
\newblock Wiley, 1986.




\bibitem{LuSqWu17}
Y.~Lu, M.~S. Squillante, C.~W. Wu.
\newblock On the control of density-dependent stochastic population processes with
  time-varying behavior.
\newblock arXiv:1709.07988, 2017.


\end{thebibliography}

%


%
%

\end{document}